\documentclass{amsart}
\usepackage[headings]{fullpage}
\usepackage{amssymb,amsmath}

\newtheorem{theorem}{Theorem}
\newtheorem{prop}[theorem]{Proposition}
\newtheorem{lemma}[theorem]{Lemma}
\newtheorem{corollary}[theorem]{Corollary}
\theoremstyle{remark}

\newtheorem*{remark}{Remark}

\newcommand{\Z}{{\bf Z}}

\newcommand{\Q}{{\bf{Q}}}
\newcommand{\ka}{{\boldsymbol k}}

\title{Nilpotency and dimension series for loops}
\author{Jacob Mostovoy}
\address{Instituto de Matem\'aticas, Unidad Cuernavaca,
 Universidad Nacional Aut\'onoma de M\'exico,
A.P. 273-3,  C.P. 62251, Cuernavaca, Morelos, MEXICO}
\email{jacob@matcuer.unam.mx}
\begin{document}

\begin{abstract}
We take a step towards the development of a nilpotency theory for
loops based on the commutator-associator filtration instead of the
lower central series. This nilpotency theory shares many essential
features with the associative case. In particular, we show that the
isolator of the $n$th commutator-associator subloop coincides with
the $n$th dimension subloop over a field of characteristic zero.
\end{abstract}

\maketitle

The lower central series for groups  can be defined in two
essentially different ways. Namely, the lower central series of a
group $G$ is a descending filtration of $G$ by normal subgroups
$$G=G_{1}\supseteq G_{2}\supseteq\ldots$$ defined inductively by
setting either: \\

$G_{i}=[G,G_{i-1}]$where $[G,H]$ is the largest of all subgroups $K$
of $G$ with the property that $K/H$ is contained in the centre of $G/H$;\\

or, $G_{i}$ to be generated by all commutators $[x,y]$ with $x\in G_{p}$
and $y\in G_{q}$ with $p+q\geq i$.\\

These two definitions are equivalent for groups. However, in the
non-associative case they give rise to rather different objects. The
first definition, with ``groups" replaced by ``loops", produces
Bruck's lower central series, see \cite{Bruck}. An analog of the
second definiton for loops was introduced in \cite{Mostovoy} under
the name of ``commutator-associator filtration". The terms of the
commutator-associator filtration contain, but do not necessarily
coincide with the corresponding terms of the lower central series.

The main advantage of the commutator-associator filtration is the
existence of a rich algebraic structure on the associated graded
abelian group, consisting of an infinite number of multilinear
operations. It can be seen that two of the operations, namely those
induced by the loop commutator and the loop associator, satisfy the
Akivis identity. However, the complete identification of this
algebraic structure is a non-trivial problem.

In this paper we set up a nilpotency theory for loops based on the
commutator-associator filtration. In this theory the  standard
techniques of the theory of nilpotent groups can be applied and
various results valid for groups can be extended to loops. In
particular, we shall prove that for an arbitrary loop the isolators
of the terms of the commutator-associator filtration coincide with
the dimension series. As a corollary, we identify the algebraic
structure  on the graded $\Q$-vector space associated to the
commutator-associator filtration: it turns out to be a Sabinin
algebra.

Throughout the text we make the emphasis on the similarities, rather
than differences, between nilpotency theories for groups and for
general loops. This should not leave the impression that extending
the nilpotency theory from groups to loops is a straightforward
task. In particular, the residual nilpotency of the free loop,
established for the lower central series by Higman \cite{Higman},
remains an open question for the commutator-associator filtration.
We did not strive for completeness; many relevant topics (such as
applications to particular classes of loops, relation to the
nilpotency of the multiplication group of the loop et cetera) have
remained outside the scope of this paper.

\subsection*{Acknowledgments}
I would like to thank Liudmila Sabinina and Jos\'e\ Mar\'\i a
P\'erez Izquierdo for discussions. This work was supported by the
CONACyT grant CO2-44100.

\section{$N$-sequences}
\subsection{The commutator-associator filtration}
The commutator-associator filtration on a loop $L$ is  defined in
terms of commutators, associators and associator deviations
\cite{Mostovoy}.

The {\em commutator} of two elements $a,b$ of $L$ is
$$[a,b]=(ba)\backslash(ab)$$
and the {\em associator} of $a,b$ and $c$ is defined by
$$(a,b,c)=(a(bc))\backslash((ab)c).$$
There is an infinite number of {\em associator deviations}. These
are functions $L^{l+3}\to L$ characterised by a non-negative number
$l$, called {\em level} of the deviation, and $l$ indices
$\alpha_{1},\ldots,\alpha_{l}$ with $0<\alpha_{i}\leq i+2$. The
deviations of level one are
$$(a,b,c,d)_{1}=((a,c,d)(b,c,d))\backslash (ab,c,d),$$
$$(a,b,c,d)_{2}=((a,b,d)(a,c,d))\backslash (a,bc,d),$$
$$(a,b,c,d)_{3}=((a,b,c)(a,b,d))\backslash (a,b,cd).$$
By definition, the deviation
$(a_{1},\ldots,a_{l+3})_{\alpha_{1},\ldots,\alpha_{l}}$ of level $l$
is equal to
$$(A(a_{\alpha_{l}})A(a_{\alpha_{l}+1}))\backslash
A(a_{\alpha_{l}}a_{\alpha_{l}+1})$$ where $A(x)$ stands for the
deviation $(a_{1},\ldots,
a_{\alpha_{l}-1},x,a_{\alpha_{l}+2},\ldots,
a_{l+3})_{\alpha_{1},\ldots,\alpha_{l-1}}$ of level $l-1$. The
associator is thought of as the associator deviation of level zero.

Now, set $\gamma_{1}L=L$ and for $n>1$ define $\gamma_{n}L$
to be the minimal normal subloop of $L$ containing
\begin{itemize}
\item{$[\gamma_{p}L,\gamma_{q}L]$ with $p+q\geq n$;}
\item{$(\gamma_{p}L,\gamma_{q}L,\gamma_{r}L)$
with $p+q+r\geq n$;}
\item{$(\gamma_{p_{1}}L,\ldots,\gamma_{p_{l+3}}L)_{
\alpha_{1},\ldots,\alpha_{l}}$ with $p_{1}+\ldots+p_{l+3}\geq n$.}
\end{itemize}
The subloop $\gamma_{n}L$ is called the {\em $n$th commutator-associator subloop} of $L$.

\begin{lemma}\label{lemma:multlin}\cite{Mostovoy}
For an arbitrary loop $L$ the commutator, the associator and the
associator deviations induce multilinear operations on the graded
abelian group $\oplus \gamma_{i}L/\gamma_{i+1}L$; these operations
respect the grading.
\end{lemma}
For the associator and the deviations the statement of the lemma
follows straight from the definition of the commutator-associator
subloops. As for the commutators, we shall now see that for
arbitrary $a\in\gamma_p L$, $b\in \gamma_q L$ and $c\in\gamma_r L$,
the commutator $[ab,c]$ is equal to $[a,c][b,c]$ modulo
$\gamma_{p+q+r}L$.

Indeed, modulo $\gamma_{p+q+r}L$
$$(ca)b\cdot [a,c]\equiv ca\cdot(b[a,c]) \equiv
ca\cdot([a,c]b)\equiv(ca)[a,c]\cdot b=(ac)b\equiv a(cb).$$ Hence,
modulo $\gamma_{p+q+r}L$
\begin{equation*}
c(ab)\cdot ([a,c][b,c]) \equiv (ca)b\cdot ([a,c][b,c])
 \equiv ((ca)b\cdot
[a,c]) [b,c]\equiv a(cb)\cdot [b,c]\equiv a\cdot(cb)
[b,c]=a(bc)\equiv (ab)c,
\end{equation*}
and therefore $([a,c][b,c])\backslash[ab,c]$ is in $\gamma_{p+q+r}L$.
Similarly one proves that  $([a,b][a,c])\backslash [a,bc]$ belongs to   $\gamma_{p+q+r}L$.
This implies that the commutator induces a bilinear operation on
$\oplus \gamma_{i}L/\gamma_{i+1}L$.

\medskip

At this point it is convenient to introduce a notion that allows to
speak of commutators, associators and deviations at the same time. A
{\em bracket of weight $n$} is an expression in $n$ indeterminates
formed by repeatedly applying  commutators, associators and
deviations, and in which every indeterminate appears only once. In
particular, the commutator is a bracket of weight 2 and a deviation
of level $l$ is a bracket of weight $l+3$. Lemma~\ref{lemma:multlin}
implies that a bracket of weight $n$ induces an $n$-linear operation
on $\oplus \gamma_{i}L/\gamma_{i+1}L$.
\begin{lemma}\label{lemma:fg}
For an arbitrary finitely generated loop $L$ the abelian groups
$\gamma_{i}L/\gamma_{i+1}L$ are finitely generated.
\end{lemma}

The group $\gamma_{i}L/\gamma_{i+1}L$ is generated by the classes of all brackets
of weight $i$. Since the brackets of weight $i$ are linear in all $i$ arguments on
$\gamma_{i}L/\gamma_{i+1}L$, the brackets of
weight $i$ whose arguments  belong to a fixed finite set of generators of $L$, are sufficient to
span $\gamma_{i}L/\gamma_{i+1}L$. However, there is only a finite number of
such brackets.

\medskip

We say that $L$ is {\em nilpotent} if there exists $n$ such that
$\gamma_{n+1}L=1$. The minimal such $n$ is called the {\em
nilpotency class} of $L$.

\begin{remark}
It can be seen that for $L$ nilpotent, the word ``normal" can be
omitted in the definition of the commutator-associator subloops.
\end{remark}

\subsection{More on deviations}
The definition of the associator deviations given above does not use
any specific property of the associator. In fact, the deviations can
be constructed for any function $\phi(x): L\to L$. We define the
{\em deviation $\phi (x_{1}, x_{2})$ derived from $\phi(x)$} by
setting
$$\phi (x_{1}, x_{2})= (\phi(x_{1})\phi(x_{2}))\backslash\phi(x_{1}x_{2}).$$
If $\phi$ is a function $L^{k}\to L$ with $k>1$ one can consider
deviations with respect to each variable. The deviations derived
from the associator are the usual associator deviations of level
one.

Now, let $w$ be a word in the free loop on $k$ generators $F_k$. For
any loop $L$ it induces a function $\phi_{w}(x_{1},\ldots, x_{k})$
from $L^{k}$ to $L$.

\begin{prop}\label{prop:lin}
Suppose that $\phi_{w}(x_{1},\ldots, x_{k})$ respects the
commutator-associator filtration on any loop, that is, $x_{i}\in
\gamma_{n_{i}}L$  implies $\phi_{w}(x_{1},\ldots, x_{k}) \in
\gamma_{n_{1}+\ldots+n_{k}}L$ for any $L$. Then a deviation derived
from $\phi_{w}$ with respect to any variable also respects the
commutator-associator filtration.
\end{prop}

The rest of this subsection is dedicated to the proof of this
statement.

Let $\phi:L^m\to L$ be a function which respects the
commutator-associator filtration. We shall say that $\phi$ is {\em
regular} if all deviations derived from $\phi$ also respect the
commutator-associator filtration. Thus, Proposition~\ref{prop:lin}
says that any function of the form $\phi_w$ that respects the
commutator-associator filtration, is regular.

Associator deviations of all levels are regular by definition. It
follows from the proof of Lemma~\ref{lemma:multlin} that the
commutator is regular.

\begin{lemma}\label{lemma:composition}
Let $\phi:L^m\to L$ be a regular function. For each $1\leq i\leq m$ let $X^i$ be a
non-empty finite set of letters and
let $\psi_i:L^{|X^i|}\to L$  be a regular function in the variables from the set $X^i$.
(We do not assume that the sets $X^i$ are disjoint.)
Then the composition $\phi(\psi_1,\ldots,\psi_m)$ is also regular.
\end{lemma}

Assume that $\phi$ is a function with two arguments and that $\psi_1$ and
$\psi_2$ are functions of the same single variable. Then for $x\in \gamma_{p}L$ and
$y\in\gamma_{q}L$ we have that, modulo $\gamma_{p+q}L$

\begin{equation*}
\begin{split}
\phi\bigl(\psi_1(xy),\psi_2(xy)\bigr)\
&=\phi\bigl(\psi_1(x)\psi_1(y)\cdot\psi_1(x,y),\psi_2(x)\psi_2(y)\cdot\psi_2(x,y)\bigr)\\
&\equiv \phi\bigl(\psi_1(x)\psi_1(y),\psi_2(x)\psi_2(y)\bigr)\\
&\equiv \Bigl(\phi\bigl(\psi_1(x),\psi_2(x)\bigr)\phi\bigl(\psi_1(x),\psi_2(y)\bigr)\Bigr)
\Bigl(\phi\bigl((\psi_1(y),\psi_2(x)\bigr)\phi\bigl((\psi_1(y),\psi_2(y)\bigr)\Bigr)\\
&\equiv \phi\bigl(\psi_1(x),\psi_2(x)\bigr)\cdot\phi\bigl(\psi_1(y),\psi_2(y)\bigr)
\end{split}
\end{equation*}
and, therefore, the deviation derived from $\phi(\psi_1,\psi_2)$
belongs to $\gamma_{p+q}L$.

The general case is entirely similar except for the complexity of
notation; we omit the proof.

\begin{lemma}\label{lemma:product}
Let $\psi_1,\psi_2: L^m\to L$ be regular functions in the same set of variables. Then the functions
$\psi_1\psi_2$, $\psi_1/\psi_2$ and  $\psi_1\backslash\psi_2$ are also regular.
\end{lemma}

For the sake of simplicity assume that $\psi_1$ and $\psi_2$ are functions in one variable.
Take $x\in \gamma_{p}L$ and $y\in\gamma_{q}L$. Then, modulo $\gamma_{p+q}L$
\begin{equation*}
\begin{split}
\psi_1(xy)\psi_2(xy)&=\bigl(\psi_1(x)\psi_1(y)\cdot\psi_1(x,y)\bigr)
\bigl(\psi_2(x)\psi_2(y)\psi_2(x,y)\bigr)\\
&\equiv \psi_1(x)\psi_1(y)\cdot\psi_2(x)\psi_2(y)\\
&\equiv \psi_1(x)\psi_2(x)\cdot\psi_1(y)\psi_2(y)
\end{split}
\end{equation*}
and it follows that $\psi_1\psi_2$ is regular.

The regularity of  $\psi_1/\psi_2$ and  $\psi_1\backslash \psi_2$ is
proved in the same manner.
The case of several variables is entirely similar.

\medskip

In order to establish the truth of Proposition~\ref{prop:lin} it is
sufficient to prove that the deviations derived from $\phi_{w}$
respect the commutator-associator filtration for nilpotent loops.
Let $L$ be of nilpotency class $N-1$. Then $\phi_w$ only depends on
the image of $w$ in $F_k/\gamma_N F_k$.

It follows from the proof of Lemma~\ref{lemma:fg} that for an arbitrary
positive integer $N$, the word $w$ can be written, modulo $\gamma_N F_k$,
as a word in  brackets of weight at least $k$ whose arguments are the generators
of $F_k$.

It can be assumed that each of these brackets contains every generator $x_1,\ldots,x_k$
of $F_k$ at least once. Indeed, the word $w$ can be written as
$$w=w_k=u_kv_k\cdot w_{k+1}$$ where $w_{k+1}$ is a word in brackets of weight
at least $k+1$, $u_k$ is a word in brackets of weight $k$ among whose arguments the generator
$x_i$ is present, and $v_k$ is a word in brackets of weight $k$ among whose arguments the
generator $x_i$ is missing. Since $\phi_{w}(x_{1},\ldots, x_{k})$ respects the commutator-associator
filtration on $F_k/\gamma_N F_k$, it follows that replacing $x_i$ by $1$ in $w$ we obtain
a word representing the identity in  $F_k/\gamma_N F_k$.
By definition, brackets of all weights respect the commutator-associator filtration. In particular,
replacing $x_i$ by $1$ in $u_k$ we also get the identity. Since $v_k$ does not change under replacing
$x_i$ by $1$, it follows that $v_k\in\gamma_{k+1}F_k/\gamma_N F_k$ and, hence, can be taken to be
equal to the identity.

This shows that $u_k$ is a word in brackets of weight $k$ which
contain all the $x_i$ among their arguments.
Lemma~\ref{lemma:composition} implies that $u_k$ gives rise to a
regular function, and therefore, by Lemma~\ref{lemma:product}
$u_k\backslash w_k$ also does. Now, writing $u_k\backslash w_k$ as
$$u_k\backslash w_k=u_{k+1}v_{k+1}\cdot w_{k+2}$$
we can repeat the argument to show that $v_{k+1}$ can be taken to be
the identity etc.

Finally, since $w=u_k(u_{k+1}(\ldots u_{N-1}))$ with all the $u_i$ giving rise to regular
functions, it follows from Lemma~\ref{lemma:product} that $\phi_w$ is also regular.

\subsection{Isolators}
Let $F$ be the free loop on the single generator $x$ and $\delta:F\to\Z$ --- the
homomorphism that sends $x$ to $1$. If $w(x)$ is a non-associative word in $x$ its
{\em degree} is defined to be the integer $\delta(w)$.

Let $L$ be a loop and $K\subset L$ -- a normal subloop. The {\em isolator} of $K$ in $L$,
denoted by $\sqrt{K}$, is the minimal normal subloop of $L$ containing all such $x\in L$
that $w(x)\in K$ for some word $w$ of non-zero degree. An element of $L$ is called
{\em periodic} if it belongs to the isolator of the identity. A loop is {\em torsion-free}
if it has no periodic elements.

\subsection{The dimension filtration}
Let us now recall the definition of the dimension filtration for loops \cite{MP}.
Let $R$ be
a commutative unital ring and $L$ --- an arbitrary loop. The {\em augmentation ideal}
$I\subset RL$ is the kernel of the $R$-linear map of the loop algebra $RL$ to $R$ that
sends every element of $L$ to $1$. The $n$th power of $I$ is the linear span of all products
of at least $n$ elements of $I$. The loop $L$ sits inside $RL$ and its intersection
with $1+I^{n}$ is a normal subloop of $L$, called the {$n$th dimension subloop over $R$}
and denoted by $D_{n}(L,R)$. In what follows we shall only consider dimension
subloops over a field $\ka$ of characteristic 0 and write $D_{n}L$ for $D_{n}(L,\ka)$.

\begin{lemma}
For any loop $L$ $\sqrt{\gamma_n L}\subseteq D_n L$.
\end{lemma}

For any loop $L$ the brackets of all weights respect the dimension
filtration. In other words, if $x\in D_{p}L$ and $y\in D_{q}L$, the
commutator $[x,y]$ belongs to $D_{p+q}L$ and similarly for the
associator and the associator deviations \cite{MP}. Since
$\gamma_{1}L=D_{1}L=L$ it follows that $\gamma_{n}L\subseteq D_{n}L$
for all $n$.

Now, let $x$ be an element of $L$ that does not belong to $D_nL$.  Then $u=1-x$
belongs to the augmentation ideal $I$ but not to $I^n$. Suppose
$u$ belongs to $I^k$ but not to $I^{k+1}$. Here $1\leq k < n$. For
any word $w$ in $x$ of degree $m\neq 0$ we have
\[w(x)=1-mu+\ldots\]
and the omitted terms are integer multiples of monomials  $u$ of degrees at least
$2$. Hence, $1-w(x)\equiv mu \mod I^{k+1}$ and, since
$\ka$ has characteristic 0, it follows that $w(x)\notin D_nL$.
However, as $\gamma_nL$ is contained in $D_nL$ for all $n$, $w(x)$
cannot be contained in $\gamma_{n}L$.

\subsection{The isolators of $\gamma_{i}L$ as an $N$-sequence}
A filtration of a loop $L$ by normal subloops $L=L_1\supseteq L_2\supseteq\ldots$
is said to be an {\em $N$-sequence} if for all $n$ any bracket of weight $n$ evaluated on
arbitrary elements $x_{i}\in L_{p_{i}}, (0<i\leq n)$ gives an element of
$L_{p_{1}+\ldots+p_{n}}$. Both the commutator-associator filtration
and the dimension filtration are $N$-sequences \cite{MP}. As in
Lemma~\ref{lemma:multlin}, the brackets of weight $n$ induce $n$-linear operations
on the graded group associated to an $N$-sequence.

\begin{prop}\label{prop:nsequence}
The filtration of any loop $L$ by $\sqrt{\gamma_n L}$ is an  $N$-sequence.
\end{prop}

The proof of Proposition~\ref{prop:nsequence} is based on the following result.

Let $\theta$ be a non-associative word on $k$ letters and let $x_{1},\ldots,x_{k}$
be elements of $L$. Define $\theta_{x_{1},\ldots, x_{k}}$ to be the subloop
of $L$ normally generated by all elements of the form
$\theta(w_{1}(x_{1}),\ldots,w_{k}(x_{k}))$ where $w_{i}$ are words on one letter.

Let $W_{i}$, where $0<i\leq k$, be words on one letter, each of non-zero degree.
\begin{lemma}\label{lemma:word}
If $L$ is nilpotent, the quotient $\theta_{x_{1},\ldots, x_{k}}/
\theta_{W_{1}(x_{1}),\ldots,W_{k}(x_{k})}$ is finite.
\end{lemma}

It is sufficient to prove Lemma~\ref{lemma:word} for the free class-$n$ nilpotent loop
$F_{k}{[n]}$ on $k$ generators $x_{1},\ldots,x_{k}$, that is, for the quotient of the
free loop on the $x_{i}$ by the $n+1$st term of its commutator-associator filtration.

The proof goes by induction on the nilpotency class. The lemma is
obvious for abelian groups. Assume it is true for free loops of
nilpotency class at most $n-1$. The kernel of the homomorphism
$F_{k}{[n]}\to F_{k}{[n-1]}$ is the commutative group
$\gamma_{n}F_{k}{[n]}$. It is enough to prove that
$\gamma_{n}F_{k}{[n]}\cap\theta_{W_{1}(x_{1}),\ldots,W_{k}(x_{k})}$
is of finite index in $\gamma_{n}F_{k}{[n]}\cap\theta_{x_{1},\ldots,
x_{k}}$.

The group  $\gamma_{n}F_{k}{[n]}$ is generated by the brackets of weight $n$
evaluated on the $x_{i}$. In particular, we can choose a basis that consists of such
brackets for the $\Q$-vector space $\gamma_{n}F_{k}{[n]}\otimes\Q$. The linearity
of the brackets implies that the homomorphism $\omega:F_{k}{[n]}\to F_{k}{[n]}$
defined by sending $x_{i}$ to $W_{i}(x_{i})$ induces a transformation of
$\gamma_{n}F_{k}{[n]}\otimes\Q$ given by a diagonal matrix with non-zero
diagonal entries. Hence,  $\omega$ induces an isomorphism of
$(\gamma_{n}F_{k}{[n]}\cap\theta_{x_{1},\ldots, x_{k}})\otimes\Q$ into itself.
On the other  hand, since the image of $\theta_{x_{1},\ldots, x_{k}}$ under
$\omega$ is contained in $\theta_{W_{1}(x_{1}),\ldots, W_{k}(x_{k})}$
it follows that $\gamma_{n}F_{k}{[n]}\cap\theta_{W_{1}(x_{1}),\ldots,
W_{k}(x_{k})}$ is of finite index in $\gamma_{n}F_{k}{[n]}\cap
\theta_{x_{1},\ldots, x_{k}}$.

\medskip

Now we are in the position to prove Proposition~\ref{prop:nsequence}.
It is sufficient to verify it for nilpotent loops; the general case
can be reduced to the case of nilpotent loops by replacing $L$ with
$L/\sqrt{\gamma_{N}L}$  with sufficiently large $N$.

Assume that $L$ is nilpotent.  Let $x\in\sqrt{\gamma_{p}L}$ and  $y\in\sqrt{\gamma_{q}L}$
so that there exist words $w_{1}$ and $w_{2}$ on one letter and of non-zero degree
such that $[w_{1}(x),w_{2}(y)]\in \gamma_{p+q}L$. Set $\theta=[x,y]$;
applying Lemma~\ref{lemma:word} we see that $\theta_{w_{1}(x),w_{2}(y)}$ is of
finite index in $\theta_{x,y}$ and, hence, there exists a word $w$ on one letter and of
non-zero degree such that $w([x,y])\in\theta_{w_{1}(x),w_{2}(y)}\subseteq \gamma_{p+q}L$.
Therefore, $[\sqrt{\gamma_{p}L},\sqrt{\gamma_{q}L}]\subseteq \sqrt{\gamma_{p+q}L}$.
Similarly one
proves that $$(\sqrt{\gamma_{p_{1}}L},\ldots,\sqrt{\gamma_{p_{l+3}}L})
_{\alpha_{1},\ldots,\alpha_{l}}\subseteq \sqrt{\gamma_{p_{1}+\ldots+p_{k+3}}L}.$$

\section{The Jennings theorem}
Now we can state our main result.
\begin{theorem}\label{theorem:jennings}
For any field $\ka$ of characteristic 0 and for any loop $L$,  the isolator
$\sqrt{\gamma_n L}$ of $\gamma_n L$ in $L$ coincides with the dimension
subloop $D_n(L,\ka)$.
\end{theorem}
The associative version of this theorem is due to Jennings \cite{Jennings}. Our
proof follows the argument given in Chapter~7 of \cite{Hall}, see also \cite{Passman}
and \cite{BirPassi}.

Theorem~\ref{theorem:jennings} implies that after tensoring with a
field of characteristic zero, the graded groups associated to the
dimension and the commutator-associator filtrations become
isomorphic. The group $\oplus D_{n}L/D_{n+1}L\otimes\Q$ has the
structure of a Sabinin algebra.

Recall that Sabinin algebras are related to Lie algebras in the same
way as loops are related to groups. They were initially introduced
by Mikheev and Sabinin as tangent structures to general affine
connections, see \cite{SM,SM-Orange}. Later, it was proved that
primitive elements in a non-associative bialgebra form a Sabinin
algebra \cite{ShU}, and that every Sabinin algebra arises this way
\cite{Perez}.

It is known from \cite{MP} that $\oplus D_{n}L/D_{n+1}L\otimes\Q$ is
the Sabinin algebra of primitive elements of the algebra $\oplus
I^{n}/I^{n+1}$, the primitive operations of Shestakov-Umirbaev
\cite{ShU} being induced by associator deviations. Therefore, we
have

\begin{corollary}
The graded group $\oplus \gamma_{n}L/\gamma_{n+1}L\otimes\Q$ is a
Sabinin algebra with the commutator and the Shestakov-Umirbaev
operations induced by the commutator and the associator deviations
on $L$ respectively.
\end{corollary}

\subsection{The outline of the proof}

We have already seen that $\sqrt{\gamma_{n}L}$ is contained $D_{n}L$. So,
just like in the associative situation, it is enough to prove that if $\sqrt{\gamma_{N}L}=1$,
then $D_{N}L$ is also trivial. We can assume that $L$ is finitely generated, since any
element of $D_{n}L$ belongs to $D_{n}L'$ where $L'\subseteq L$ is some finitely
generated subloop.

\medskip

Let us fix some notation. For $a\in L$ denote the corresponding element of the left
multiplication group of $L$ by $\lambda_{a}$. Similarly, for any $v$ in the loop
algebra $\ka L$ we write $\lambda_{v}$ for the
left multiplication by $v$ in $\ka L$. Writing a product
$a_{1}a_{2}\ldots a_{m}$ without parentheses we mean
$a_{1}(a_{2}(\ldots a_{m}))=\lambda_{a_{1}}\lambda_{a_{2}}
\ldots\lambda_{a_{m}}(1)$. The expression $a^{m}$ will stand for
$\lambda_{a}^{m}(1)$.

\medskip

All the quotients $\sqrt{\gamma_{n}L}/\sqrt{\gamma_{n+1}L}$ are
torsion-free; let $M$ be the sum of the ranks of these quotients. There are $x_{i}\in L$
with $1\leq i\leq M$ and integers $c_{j}$ with $1\leq j \leq N$ such that
$c_{1}=1$, $c_{j}\leq c_{j+1}$, $c_{N}=M+1$ and
$$\sqrt{\gamma_{n}L}=
\{\sqrt{\gamma_{n+1}L}, x_{c_{n}}, x_{c_{n}+1},\ldots,x_{c_{n+1}-1}\}.$$
Then each element of $L$ can be uniquely written as
$\lambda_{x_{1}}^{r_{1}}\lambda_{x_{2}}^{r_{2}}\ldots
\lambda_{x_{M}}^{r_{M}}(1)$ with $r_{i}$ integers.

Let $x$ be a generator of the infinite cyclic group and set $u=1-x$.
By Lemma~7.2 of \cite{Hall}, the group ring of the infinite cyclic group has the basis
consisting of $1,u,u^{2},\ldots$ together with $u^{m}x^{-1}, u^{m}x^{-2},\ldots$
where $m$ is any positive integer.  Let $u_{i}=1-x_{i}\in \ka L$, then
$\lambda_{u_{i}}=1-\lambda_{x_{i}}$. Therefore, we have the following
\begin{lemma}\label{lemma:basis}
The loop algebra $\ka L$ has a basis consisting of all elements of the form
$$\Lambda_{1}\Lambda_{2}\ldots\Lambda_{M}(1)$$ where $\Lambda_{i}$ is equal to
either $\lambda_{u_{i}}^{r_{i}}$ or $\lambda_{x_{i}}^{-s_{i}}\lambda_{u_{i}}^{N}$
with $r_{i}$ a non-negative and $s_{i}$ --- a positive integer.
\end{lemma}

Define $\mu(u_{i})$ to be the largest number $k$ such that
$x_{i}\in\sqrt{\gamma_{k}L}$. For any basis element $v$ of the form described in
Lemma~\ref{lemma:basis} we define $\mu(v)\geq N$ if at least one of the $\Lambda_{i}$
in $v$ has the form  $\lambda_{x_{i}}^{-s}\lambda_{u_{i}-1}^{N}$; otherwise
$v$ is of the form $v= \lambda_{u_{1}}^{r_{1}}\ldots\lambda_{u_{M}}^{r_{M}}$
and we set $\mu(v)=\sum \mu(u_{i})r_{i}$.

For $k\leq N$ denote by $E_{k}$ the vector space over $\ka$ spanned by those basis
elements $v$ with $\mu(v)\geq k$; for $k>N$ set $E_{k}=E_{N}$.  It is clear that
$E_{k}$ is contained in the $k$th power $I^{k}$ of the augmentation ideal $I$.

\begin{lemma}\label{lemma:ekl}
$E_{k}E_{l}\subseteq E_{k+l}.$
\end{lemma}

The proof of this lemma will be given in the next section. Now,
assuming the the truth of Lemma~\ref{lemma:ekl}, let us finish the
proof of Theorem~\ref{theorem:jennings}.

The augmentation ideal $I$ is the same thing as $E_{1}$. Since $E_{k}\subseteq I^{k}$,
it follows from Lemma~\ref{lemma:ekl} that $I^{k}$  coincides with $E_{k}$.
Lemma~\ref{lemma:basis} implies that the elements $u_{i}=\lambda_{u_{i}}(1)$
with $\mu(u_{i})=k$ (these are the $u_{i}$ with $c_{k}\leq i< c_{k+1}$) are linearly
independent modulo $E_{k+1}$, and hence, modulo $I^{k+1}$.

Now, any $y\in\sqrt{\gamma_{k}L}-\sqrt{\gamma_{k+1}L}$  is of the form
$\lambda_{x_{j}}^{r_{0}}\lambda_{x_{j+1}}^{r_{1}}\ldots
\lambda_{x_{j+l}}^{r_{l}}(z)$ with $j=c_{k}$,  $l=c_{k+1}-c_{k}-1$,
$z\in\sqrt{\gamma_{k+1}L}$ and not all $r_{i}$ equal to zero. Since $1-z$ is in
$I^{k+1}$, we have
\[1-y\equiv \sum r_{i}u_{j+i}\ \mod I^{k+1}\]
where the sum is over all $0\leq i\leq l$. Therefore, $1-y$ does not belong to $I^{k+1}$
and, hence, $y\notin D_{k+1}L$. In particular, if $y\neq 1$, then $y\notin D_{N}L$.

\section{Proof of Lemma~\ref{lemma:ekl}}

Given a set $X$ of elements of $L$, an {\em elementary bracket} with
respect to $X$ is a bracket whose arguments belong to $X$. Take
$X=\{a, b_{1},\ldots,b_{p}, c_{1},\ldots, c_{q}\}$. Then, using the
definition of associator deviations, one can decompose the loop
associator $(a, b_{1} \ldots b_{p}, c_{1}\ldots c_{q})$ as a product
of elementary brackets. For every pair of non-empty subsets
$I=\{i_{1},\ldots i_{|I|}\}\subseteq\{1,2,\ldots,p\}$ and
$J=\{j_{1},\ldots j_{|J|}\}\subseteq\{1,2,\ldots,q\}$ this product
contains precisely one deviation of the form $(a,
b_{i_{1}},\ldots,b_{i_{|I|}},c_{j_{1}},\ldots,c_{j_{|J|}})_{*,\ldots,*}$.
We shall fix once and for all such a decomposition and call it
$\Pi$.

Let $S$ be some subset of the set of all elementary brackets that
form the product $\Pi$. One can then form a product of elementary
brackets $\Pi_{S}$ by deleting from $\Pi$ all the brackets that do
not belong to $S$. Now, replace in the product $\Pi_{S}$ each
elementary bracket $w$ by $w-1$; the resulting element of $\ka L$ is
denoted by  $P_{S}$. If  $S$ is empty, then $\Pi_{S}=1$ and
$P_{S}=0$.

Write $A, B_i$ and $C_j$ for $1-a$,  $1-b_i$ and  $1-c_j$
respectively. Denote by $B_I$ and $C_J$ the products
$B_{i_1}B_{i_2}\ldots B_{i_{|I|}}$ and  $C_{j_1}C_{j_2}\ldots
C_{j_{|J|}}$ respectively, where  $I=\{i_{1},\ldots i_{|I|}\}$ and
$J=\{j_{1},\ldots j_{|J|}\}$, in a similar way we define products
$b_I$ and $c_J$. If $I$ (or $J$) is empty, then $B_I=1$ ($C_J=1$,
respectively).

Then the following formula holds:
\begin{multline}\label{formula:ass}
A\cdot(B_{1}B_{2}\ldots B_{p}\cdot C_{1} C_{2}\ldots C_{q}) -\bigl(A
B_{1}B_{2}\ldots B_{p}\bigr)\cdot
C_{1} C_{2}\ldots C_{q}\\
= (-1)^{p+q}\sum_{S}\Bigl( \sum_{I,J} (-1)^{|I|+|J|} a \cdot B_{I}
C_{J} \Bigr) \cdot P_{S}.
\end{multline}
Here the sum inside the brackets on the right-hand side is taken
over all subsets $I\subseteq\{1,\ldots, p\}$ and
$J\subseteq\{1,\ldots, q\}$ with the property that if no bracket in
$S$ contains $b_i$ (or $c_j$) as an argument, then $i\in I$ (or
$j\in J$, respectively).

In order to prove (\ref{formula:ass}), notice that
$$a b_I\cdot c_J=(a\cdot b_I c_J)\bigl(1+\sum_{S\subseteq S_{I,J}} P_S\bigr),$$
where $S_{I,J}$ is the subset consisting of all brackets from $\Pi$
which contain only the variables $a$, $b_I$ and $c_J$. Also,
\begin{equation*}
\begin{split}
a B_I\cdot C_J\ &=\sum_{I'\subseteq I,J'\subseteq J}
(-1)^{|I'|+|J'|} a b_{I'}\cdot c_{J'}\\
&= \sum_{I'\subseteq I,J'\subseteq J} (-1)^{|I'|+|J'|} (a \cdot
b_{I'}c_{J'})\Bigl(1+\sum_{S\subseteq S_{I',J'}}
P_S\Bigr) \\
&=a\cdot B_I C_J + \sum_{I'\subseteq I,J'\subseteq J}
(-1)^{|I'|+|J'|} (a\cdot b_{I'}c_{J'}) \sum_{S\subseteq S_{I',J'}}
P_S.
\end{split}
\end{equation*}
Now, $a B_I\cdot C_J-a\cdot B_I C_J=A\cdot B_I C_J-A B_I\cdot C_J$.
It remains to calculate the coefficient at $P_S$ for given $S$:
\begin{equation*}
\sum_{I'\subseteq I''\subseteq I, J'\subseteq J''\subseteq
J}(-1)^{|I''|+|J''|} (a\cdot b_{I''}c_{J''})=\sum_{I-I''\subseteq
I', J-J''\subseteq J'} (-1)^{|I-I''|+|J-J''|}(a\cdot
B_{I''}C_{J''}).
\end{equation*}
Now, setting $I=\{1,\ldots,p\}$ and $J=\{1,\ldots,q\}$  and writing
$I,J$ instead of $I'',J''$ we get (\ref{formula:ass}).

\medskip

We shall need two other formulae similar to (\ref{formula:ass}).
Consider the {\em anti-associator}
$$(a,b,c)'=((ab)c)\backslash (a(bc)).$$
Mimicking the construction of deviations for the associator, we can
build the hierarchy of deviations of all levels derived from the
anti-associator. Then we have the following formula:

\begin{multline}\label{formula:anti-ass}
\bigl(A_{1}A_{2}\ldots A_{p} \cdot B \bigr)\cdot
C_{1} C_{2}\ldots C_{q}-
A_{1}A_{2}\ldots A_{p}\cdot B C_{1} C_{2}\ldots C_{q}\\
= (-1)^{p+q}\sum_{S}\Bigl( \sum_{I,J} (-1)^{|I|+|J|} A_{I}
 b \cdot C_{J} \Bigr)\cdot Q_{S}.
\end{multline}

Here $Q_S$ is defined exactly as $P_S$ but with anti-associators and
the deviations derived from them instead of associators and
associator deviations. All other symbols have the same meaning as in
$(\ref{formula:ass})$. In the particular case when $p=1$, the
formulae $(\ref{formula:anti-ass})$ and $(\ref{formula:ass})$ give
\begin{align}\label{formula:anti-ass2}
A B \cdot C_{1} C_{2}\ldots C_{q}- A\cdot B C_{1} C_{2}\ldots C_{q}
&= (-1)^{q+1}\sum_{S}\Bigl( \sum_{J} (-1)^{|J|} a b \cdot C_{J}
\Bigr)\cdot Q_{S}\\
\label{formula:ass2} &= (-1)^{q}\sum_{S}\Bigl( \sum_{J} (-1)^{|J|}
a\cdot bC_{J} \Bigr)\cdot P_{S}
\end{align}
where $A$ and $B$ stand for $A_1$ and $B_1$ respectively, and $a$
and $b$ --- for $a_1$ and $b_1$.

The same construction can also be performed for the commutator. The
resulting formula reads
\begin{equation}\label{formula:comm}
A_{1}A_{2}\ldots A_{p}\cdot B -B\cdot A_{1}A_{2}\ldots A_{p} =
(-1)^{p+1}\sum_{S}\Bigl( \sum_{I} (-1)^{|I|}b \cdot A_{I} \Bigr)
\cdot R_{S}.
\end{equation}

\begin{remark}
The anti-associator is readily seen to respect the
commutator-associator filtration. Proposition~\ref{prop:lin} then
implies that the deviations of all levels derived from the
anti-associator also respect the commutator-associator filtration.
The same thing can be said about the commutator and the deviations
of all levels derived from it.
\end{remark}

\medskip

Set $v_i=1-x_i\backslash 1$ and consider products of the form
\begin{equation}\label{form}
\Lambda'_{1}\Lambda'_{2}\ldots\Lambda'_{M}(1)
\end{equation}
where $\Lambda'_{i}$ is equal to
$\lambda_{v_{i}}^{d_{i}}\lambda_{u_{i}}^{r_i}$ with $d_{i},
r_{i}\geq 0$. For every such product $u$ set
$\mu'(u)=\sum{\mu(u_i)r_i}$ and $l(u)=\sum{(d_i+r_i)}$.

Let $E'_{k}$ be the subspace of $\ka L$ spanned by the products
$\Lambda'_{1}\ldots\Lambda'_{M}(1)$ with $\mu'\geq k$. For $1\leq
s\leq M$, let $E'_{k, s}$  be the subspace of $E'_{k}$ spanned by
the products that have
$\Lambda'_{1}=\Lambda'_{2}=\ldots=\Lambda'_{s-1}=1$. In particular,
$E'_{k,1}=E'_k$.

\begin{lemma}\label{lemma:key}
If $y\in E'_{k,s}$, $z\in E'_{l,s}$ and
$x=\lambda^{d}_{v_{s-1}}\lambda^{r}_{u_{s-1}}(y)$, then
\begin{equation}\label{conditions}
x z=\sum_i \alpha_i
\lambda^{\delta_i}_{v_{s-1}}\lambda^{\rho_i}_{u_{s-1}}(w_i)
\end{equation}
where $w_i\in E'_{m_i,s}$, $\alpha_i$ are integers,
$m_i=k+l+\mu(u_{s-1})(r-r_i)$, $\delta_i\leq d$ and $\rho_i \leq r$.

In particular, for  $y\in E'_{k,s}$ and $z\in E'_{l,s}$  the product
$yz$ is contained in $E'_{k+l,s}$.
\end{lemma}

\begin{proof}
Let $z\in E'_{l,t}$ where $t\geq s$; we shall use descending
induction on $t$. For $t=M$ the statement of the lemma is obvious
since all elements of $E'_{l,M}$ commute and associate with
everything. Suppose that for $t>q$ and all $y,z$ as above the lemma
has been established.

Consider first the case $z=u_q$. Set $p=l(x)$ and use induction on
$p$.

For $p=1$ we have that $x$ is equal to either $u_j$ or $v_j$ for
some $j\geq s-1$, and there are two possibilities. If $j\leq q$, the
condition (\ref{conditions}) is satisfied automatically. If $j>q$ we
have
$$u_{j}u_q-u_qu_{j}=\bigl((1-u_q)(1-u_{j})\bigr)\bigl([x_q,x_{j}]-1\bigr)$$
and
$$v_{j}u_q-u_qv_{j}=\bigl((1-u_q)(1-v_{j})\bigr)\bigl([x_q,x_{j}\backslash 1]-1\bigr).$$
(These formulae are particular cases of (\ref{formula:comm}).) Since
$[x_q,x_{j}]-1$ $[x_q,x_{j}\backslash 1]-1$ belong to
$E'_{\mu'(u_q)+\mu'(u_j), q+j}$ the first induction assumption
implies (\ref{conditions}) for $p=1$.

Assume that (\ref{conditions}) holds for $z=u_q$ and all $p<p_0$.

Take $y\in E'_{k,s}$ with
$l(\lambda_{v_{s-1}}^{d}\lambda_{u_{s-1}}^{r}(y))=p_0$. If $d>0$,
write $y=v_{s-1}\tilde{y}$. Then $v_{s}\tilde{y}\cdot
u_{q}-v_{s}\cdot \tilde{y}u_{q}$ can be re-written  with the help of
(\ref{formula:ass}). The right-hand side of (\ref{formula:ass}) has
the form $$\pm\sum_{S,I} (-1)^{|I|} \Bigl((x_s\backslash 1) \cdot
\tilde{y}_{I}x_q \Bigr) \cdot P_{S}.$$ By the second induction
assumption, the product $\tilde{y}_{I} u_q$ satisfies
(\ref{conditions}), and $P_S\in E'_{*, q'}$ with $q'>q$. Therefore,
applying the first induction assumption we see that the product
$v_{s}\tilde{y}\cdot u_{q}$ also satisfies (\ref{conditions}), so we
see that (\ref{conditions}) holds for $z=u_q$ and $p=p_0$ as well.
In the situation where $d=0$ but $r\neq 0$ the argument is
completely analogous.

If $d=0$ and $r=0$, we only have to consider the case when $q=s$.
Applying (\ref{formula:comm}) and the induction assumptions we see
again that (\ref{conditions}) is fulfilled, so the lemma holds
whenever $z=u_q$.

In a similar fashion one verifies the lemma for $z=v_q$.

Let us now pass to the case of arbitrary $z\in E'_{l,q}$. If the
condition (\ref{conditions}) fails for some
$z=\lambda_{u_q}^{c}(z')$ (where $z\in E'_{l-c\mu'(u_q),q+1}$) and
some $d,r$ and $y$, choose the counterexample with the smallest
possible $c$. Then, on one hand, $xu_q\cdot\lambda_{u_q}^{c-1}(z')$
satisfies (\ref{conditions}). On the other hand,
$xu_q\cdot\lambda_{u_q}^{c-1}(z')-x\cdot\lambda_{u_q}^{c}(z')$ can
be re-written using (\ref{formula:anti-ass}). However, using the
induction assumption, we see that the right-hand side of
(\ref{formula:anti-ass}) is a linear combination of products of the
form $\lambda^{\delta_i}_{v_{s-1}}\lambda^{\rho_i}_{u_{s-1}}(w_i)$
with $\delta_i\leq d$ and $\rho_i \leq r$ and the $w_i$ belonging to
the ``correct'' terms of the filtration $E'_{*,*}$. Therefore, no
such counterexample can exist.

Finally, it may happen that (\ref{conditions}) fails for some
$z=\lambda_{v_q}^{b}\lambda_{u_q}^{c}(z')$. Then the argument of the
previous paragraph carries over to this situation without
modifications. This completes the induction step.

\end{proof}

\begin{lemma}\label{lemma:wellwellwell}
If $u\in E'_{k,s}$ and $v\in E'_{l,s}$, then
$$(x_s u) v=\sum_i \alpha_i x_s w_i$$ where $w_i\in E'_{k+l,s}$ and $\alpha_i$ are
integers.
\end{lemma}

\begin{proof}
Assume that $v\in E'_{l,t}$ and use descending induction on $t$. For
$t=M$ there is nothing to prove since $x_M$ commutes and associates
with everything. Suppose that the lemma is established for $t>q$.

We have $$(x_s u)v-x_s (uv)=u_s (uv)-(u_s u)v$$ which, by
formula~(\ref{formula:ass}) is equal to $\pm \sum_{S}\Bigl(
\sum_{I,J} (-1)^{|I|+|J|} x_s \cdot u_{I} v_{J} \Bigr) \cdot P_{S}$.
It follows from Lemma~\ref{lemma:key} that $uv\in E'_{k+l,s}$.
Moreover, in each term $u_{I} v_{J}\in E'_{k',s}$ and $P_S\in
E'_{k'',s'}$ with $k'+k''\geq k+l$ and $s'>q$ and, hence, by the
induction assumption $u_s (uv)-(u_s u)v\in E'_{k+l,s}$.

\end{proof}

For $1\leq s\leq M$, let $E_{k, s}$  be the subspace of $E_{k}$
spanned by the basis elements
$\Lambda_{1}\Lambda_{2}\ldots\Lambda_{M}(1)$ that have
$\Lambda_{1}=\Lambda_{2}=\ldots=\Lambda_{s-1}=1$. We have
\begin{lemma}\label{lemma:final}
$E_{k,s}=E'_{k,s}$.
\end{lemma}
Together with Lemma~\ref{lemma:key} this establishes
Lemma~\ref{lemma:ekl} since $E_k=E_{k,1}$.

\begin{proof}
Let us prove that every product of the form
$\Lambda'_{s}\Lambda'_{s+1}\ldots\Lambda'_{M}(1)$ with $\mu'\geq k$
is in $E_{k,s}$. The proof uses descending induction on $s$. For
$s=M$ there is nothing to prove.

Assume that this statement is true for $s=q+1$. If it does not hold
for $s=q$, the set of all $d$ such that
$\lambda_{v_{q}}^{d}\lambda_{u_{q}}^{r}(v)$ does not belong to
$E_{k, q}$ for some $r\geq 0$ and some $v\in E_{k-r\mu(u_q),q+1}$,
is non-empty. Take the smallest such $d$; clearly, $d\geq 1$. Choose
$w=\lambda_{v_{q}}^{d-1}\lambda_{u_{q}}^{r}(v)\in E_{k,q}$ such that
$v_q w$ is not contained in $E_{k,q}$. We have
\begin{equation}\label{formula:simple}
w-x_q((x_q\backslash 1)w)=(x_q\cdot x_q\backslash 1)w
-x_q((x_q\backslash 1)w)=u_q v_{q} \cdot w - u_q\cdot v_{q} w.
\end{equation}
Since $w$ can be written as a product $w_1\ldots w_m$ with $w_i$ of
the form $u_j$ or $v_j$, the formula~(\ref{formula:anti-ass2}) can
be applied to the expression on the right-hand side of
(\ref{formula:simple}); it is equal to
\begin{equation*}
\pm \sum_{S}\Bigl( \sum_{J} (-1)^{|J|} w_{J} \Bigr)\cdot Q_{S}.
\end{equation*}
This lies in $E_{k,q}$ by Lemma~\ref{lemma:key} and the induction
assumption; so does $w$. Therefore $x_q((x_q\backslash1)w)$ also is
in $E_{k,q}$ and, hence, $x_q(v_q w)\in E_{k,q}$.

The operator $\lambda^{-1}_{x_q}$ of left division by $x_q$ can be
written as a linear combination of operators of the form
$\lambda^{t}_{u_q}$ and $\lambda^{-1}_{x_q}\lambda^{N}_{u_q}$ with
$t\geq 0$, so $\lambda_{x_q}^{-1}(x_q(v_q w))=v_q w$ is also in
$E_{k,q}$, which contradicts our choice of $w$. Therefore, for $s=q$
all products of the form
$\Lambda'_{s}\Lambda'_{s+1}\ldots\Lambda'_{M}(1)$ with $\mu'\geq k$
belong to $E_{k,q}$.

\medskip

It remains to see that $E_{k,s}\subseteq E'_{k,s}$. Assume that we
have this established for $s>q$. Consider the set of all $d$ such
that $\lambda_{x_{q}}^{-d}\lambda_{u_{q}}^{r}(v)$ does not belong to
$E'_{k, q}$ for some $v\in E_{k',q+1}$ and some $r\geq 0$ with
$k'+r\geq k$. If this set is empty we are done. If not, take the
smallest such $d$; we have $d\geq 1$. Choose
$w=\lambda_{x_{q}}^{-d+1}\lambda_{u_{q}}^{r}(v)\in E_{k,q}$ such
that $x_q\backslash w$ is not contained in $E'_{k,q}$.

By construction $w\in E'_{k,q}$. Applying (\ref{formula:ass2}) to
the expression on the right-hand side of (\ref{formula:simple})
gives
$$u_q v_{q} \cdot w - u_q\cdot v_{q} w=
\pm \sum_{S}\sum_{J} (-1)^{|J|}\Bigl( x_q w_{J} - x_q\cdot v_q
w_{J}\Bigr)\cdot P_{S}$$ It follows then from
Lemma~\ref{lemma:wellwellwell} that
\begin{equation*}
w-x_q((x_q\backslash 1)w)=\sum_i \alpha_i x_q \tilde{w}_i,
\end{equation*}
where each $\tilde{w}_i$ is in $E'_{k,q}$ and the $\alpha_i$ are
some integers. Therefore, $x_q\backslash w$ is also in $E'_{k,q}$.
\end{proof}

\end{document}